\input amssym.def 
\input amssym
\magnification=1200
\parindent0pt
\hsize=16 true cm
\baselineskip=13  pt plus .2pt
$ $

\def\Z{{\Bbb Z}}
\def\D{{\Bbb D}}
\def\A{{\Bbb A}}
\def\S{{\Bbb S}}
\def\C{{\cal C}}
\def\R{{\Bbb R}}
\def\H{{\cal H}}

\centerline {\bf On large orientation-reversing finite group-actions on 3-manifolds}

\centerline {\bf and equivariant Heegaard decompositions}

\bigskip \bigskip

\centerline {Bruno P. Zimmermann}

\medskip

\centerline {Universit\`a degli Studi di Trieste}

\centerline {Dipartimento di Matematica e Geoscienze}

\centerline {34127 Trieste, Italy}

\bigskip \bigskip

{\bf Abstract.}  We consider finite group-actions on closed, orientable and
nonorientable 3-manifolds; such a finite group-action leaves invariant the
two handlebodies of a Heegaard splitting of $M$ of some genus $g$. 
The maximal possible order of a finite group-action of an orientable or
nonorientable handlebody of genus $g>1$ is
$24(g-1)$, and in the present paper we characterize the 3-manifolds $M$ and
groups $G$ for which the maximal possible order $|G| = 24(g-1)$ is obtained,
for some $G$-invariant Heegaard splitting of genus $g>1$.  If
$M$ is reducible then it is obtained by doubling an action of maximal
possible order
$24(g-1)$ on a handlebody of genus $g$. If $M$ is irreducible then it is a
spherical, Euclidean or hyperbolic manifold obtained as a quotient of one of
the three geometries by a normal subgroup of finite index of a Coxeter group
associated to a Coxeter tetrahedron, or of a twisted version of such a
Coxeter group.

\bigskip \bigskip

{\bf 1. Introduction}

\medskip

The maximal possible order of a finite group $G$ of
orientation-preserving diffeomorphisms of an orientable handlebody of
genus $g>1$ is $12(g-1)$ ([Z1], [MMZ]); for orientation-reversing finite
group-actions on an orientable handlebody and for actions on a
nonorientable handlebody, the maximal possible order is
$24(g-1)$; we will always assume $g>1$ in the present paper.

\medskip

Let $G$ be a finite group of diffeomorphisms of a closed, orientable or
nonorientable 3-manifold $M$. We define the (equivariant) Heegaard genus of
such a $G$-action as the minimal genus $g>1$ of a Heegaard decomposition of
$M$ into two handlebodies of genus $g$ (nonorientable if $M$ is
nonorientable) such that both handlebodies
are invariant under the $G$-action. Then  $|G| \le 24(g-1)$, and in the
maximal case $|G| = 24(g-1)$ we call both the $G$-action and the 3-manifold
{\it strongly maximally symmetric} (the term maximally symmetric is used
in various papers for the case of orientation-preserving
actions of maximal possible order $12(g-1)$ on orientable 3-manifolds, see
[Z2], [Z7] or the survey [Z4]).  In the present paper we  characterize the strongly
maximally symmetric finite group-actions, using the approach to finite
orientation-reversing group-actions on handlebodies in [Z3].

\medskip

In order to state our results, we introduce some notation (see [T1], [T2]
for the following). A {\it Coxeter tetrahedron} is a tetrahedron all of whose
dihedral angles are of the form
$\pi/n$ (denoted by a label $n$ of the edge, for some integer $n \ge 2$) and,
moreover, such that at each of the four vertices the three angles of the 
adjacent edges define a spherical triangle (i.e., $1/n_1 + 1/n_2 + 1/n_3 >
1$).  Such a Coxeter tetrahedron can be realized as a spherical, Euclidean or
hyperbolic tetrahedron in the 3-sphere $S^3$, Euclidean 3-space $\R^3$ or
hyperbolic 3-space $\Bbb H^3$, and will be denoted by $\C(n,m;a,b;c,d)$ where
$(n,m), (a,b)$ and $(c,d)$ are the labels of pairs of opposite edges. We
denote by
$C(n,m;a,b;c,d)$ the {\it Coxeter group} generated by the
reflections in the faces of $\C(n,m;a,b;c,d)$, a properly discontinuous
group of isometries of one of these three geometries.

\medskip

We are interested in particular in the Coxeter tetrahedra 
$\C(n,m) = \C(n,m;2,2;2,3)$ and the corresponding Coxeter groups 
$C(n,m) = \C(n,m;2,2;2,3)$ which are exactly the following:

\medskip

spherical:  \hskip 5mm $C(2,2), \;\;  C(2,3), \;\;  C(2,4), \;\; C(2,5),
\;\;  C(3,3), \;\;  C(3,4), \;\;  C(3,5)$;
 
Euclidean:    \hskip 3mm $C(4,4)$;

hyperbolic:  \hskip 2mm $C(4,5), \;\;   C(5,5).$

\medskip

We consider also the Coxeter tetrahedra of type 
$\C(n,m;3,3;2,2)$; such a Coxeter tetrahedron has a rotational symmetry
$\tau$ of order two which exchanges the opposite edges with labels 3 and 
2 and acts as an inversion on the two edges with lables $n$ and $m$. The
involution $\tau$ can be realized by an isometry and hence defines a group of
isometries $C_\tau(n,m)$ containing the Coxeter group 
$C(n,m;3,3;2,2)$ as a subgroup of index two; we call $C_\tau(n,m)$
a {\it twisted Coxeter group}. The twisted Coxeter groups of type $C_\tau(n,m)$
are the following:

\medskip

spherical:  \hskip 5mm  $C_\tau(2,2), \;\;  C_\tau(2,3), \;\; C_\tau(2,4)$;

Euclidean:    \hskip 3mm $C_\tau(3,3)$;

hyperbolic:  \hskip 2mm $C_\tau(2,5), \;\; C_\tau(3,4), \;\;
C_\tau(3,5), \;\;  C_\tau(4,4),  \;\; C_\tau(4,5), \;\;  C_\tau(5,5).$

\bigskip

Our main result is:

\bigskip

{\bf Theorem.}   {\sl  i) A reducible, strongly maximally symmetric
$G$-manifold $M$ is obtained by doubling a $G$-action of maximal possible
order $24(g-1)$ on an orientable or nonorientable handlebody of genus
$g>1$ (i.e., by taking the double along the boundary of both the handlebody
and its $G$-action).

\medskip

ii) An irreducible, strongly maximally symmetric
$G$-manifold $M$ is spherical, Euclidean or hyperbolic and  obtained as a
quotient of the 3-sphere, Euclidean or hyperbolic 3-space by a normal
subgroup of finite index, acting freely, of a Coxeter group $C(n,m)$ or a
twisted Coxeter group $C_\tau(n,m)$; the
$G$-action is obtained as the projection of the Coxeter or twisted
Coxeter group to $M$.}

\bigskip

There are only finitely many possibilities in the
spherical case, and most of the finite group-actions are on the
3-sphere. The genera  of the strongly maximally symmetric group-actions on
$S^3$ can be computed from the orders of the spherical Coxeter groups (see
[CM], [Z2, Table 1]), and one obtains:

\bigskip

{\bf Corollary.}   {\sl The genera of the strongly maximally symmetric
finite group actions on the 3-sphere are $g$ = 2, 3, 5, 11, 6, 17 and
601 in the untwisted cases, and $g$ = 4, 11, and 97 in the twisted
cases.}

\bigskip

For orientation-preserving actions of maximal possible order $12(g-1)$ on
the  3-sphere, the genera are determined in [WZ1, Theorem 3.1] and [WZ2], and
there are in addition the values $g$ = 9, 25, 121 and 241.

\medskip

In section 3 we consider also the case of $G$-actions
of maximal possible order $48(g-1)$, allowing $G$-actions which
interchange the two handlebodies of a Heegaard splitting of a 3-manifold $M$.
We use methods from [Z3], in particular we correct and complete 
results in [Z3] where the twisted cases are missing. Our results for the
spherical case are obtained also in [WWZ] where more general group actions of large
orders on the 3-sphere are considered.

\bigskip

{\bf 2. Proof of the Theorem}

\medskip

Let $G$ be a finite group of maximal possible order $24(g-1)$
which acts on a closed 3-manifold $M$ and leaves invariant the handlebodies
$V_g$ and
$V_g'$ of a Heegaard decomposition $M = V_g \cup_{\partial}V_g'$ of genus
$g>1$.  By [Z3],  each of $V_g/G$ and $V_g'/G$ is a
{\it handlebody orbifold} obtained by glueing two 3-disk orbifolds $D^3/G_1$ and
$D^3/D_2$ (quotients of the closed 3-disk $D^3$ by spherical groups $G_1$ and
$G_2$) along a common 2-disk suborbifold
$D^2/\D_n$ of their boundaries (a quotient of the 2-disk $D^2$ by a dihedral
group $\D_n$ of order $2n$), with orbifold Euler characteristic 
$\;\; 1/|\D_n|  -  1/|G_1| -  1/|G_2|$  (see [T1], [T2] for basic facts about
orbifolds, and [Z7] for the orientation-preserving case).  For the case of maximal
possible order
$24(g-1)$, there are exactly eight such handlebody orbifolds 
which are the handlebody orbifolds of orbifold Euler characteristic -1/24, 
the largest possible value smaller than 0; since the orbifold Euler
characteristic is multiplicative under finite orbifold coverings, $(-1/24)|G|
= 1-g$, $|G| = 24(g-1)$. These minimal handlebody orbifolds can best be codified by
their orbifold fundamental groups
$G_1*_{\D_n}G_2$ which, by [Z3, Theorem 1], are exactly the following eight
free products with amalgamation:

$$\bar\D_2*_{\bar\Z_2}\bar\D_3, \;\;\; \bar\D_3*_{\bar\Z_3}{\bar\A_4}, \;\;\;
\bar\D_4*_{\bar\Z_4}{\bar\S_4}, \;\;\;  \bar\D_5*_{\bar\Z_5}{\bar\A_5},$$ 
$$\D_{2*2}*_{\bar\Z_2}\bar\D_3, \;\;\; \D_{2*3}*_{\bar\Z_3}{\A_4},
\;\;\;
\D_{2*4}*_{\bar\Z_4}{\bar\S_4}, \;\;\; 
\D_{2*5}*_{\bar\Z_5}{\bar\A_5};$$

\medskip

here  $\bar \D_n = [2,2,n]$, 
$\bar\A_4 =[2,3,3]$, $\bar\S_4 = [2,3,4]$ and $\bar\A_5 =
[2,3,5]$ denote the extended dihedral, tetrahedral,
octahedral and dodecahedral groups (generated by the reflections in the
corresponding spherical triangles with angles $\pi/m$).  The group
$\D_{2*n}$ is a spherical group of order $4n$ (isomorphic to the dihedral
group $\D_{2n}$), a subgroup of index two in the extended dihedral group $\bar
\D_{2n}$, with the standard dihedral action by rotations and
reflections on the equatorial section $S^1$ of the 2-sphere $S^2$, the 
reflections corresponding alternatingly to rotations  and  reflections in
great circles of $S^2$.  The extended cyclic group $\bar Z_n$ of order $2n$ 
is isomorphic to the dihedral group $\bar D_n$ and  acts in the standard way
by rotations and reflections on $D^2$, and also on $D^3$.

\medskip

For example, in  the case of $\bar\D_5*_{\bar\Z_5}{\bar\A_5}$ the associated
handlebody orbifold $\H(\bar\D_5*_{\bar\Z_5}{\bar\A_5})$ is obtained as
follows. The 3-disk orbifold $D^3/\bar\D_5$ is a cone over its
orbifold boundary,  the spherical 2-orbifold $S^2/\bar\D_5$ which is just a
triangle with angles $\pi/2, \pi/2$ and $\pi/5$; the singular points of this
triangle consists of its three sides which are reflection axes with local
groups $\Z_2$, seperated by the three vertices with local groups $\D_2$,
$\D_2$ and  $\D_5$. In a similar way, the 3-disk orbifold
$D^3/\bar \A_5$ is constructed. In their boundaries, both 3-disk orbifolds
$D^3/\bar\D_5$ and
$D^3/\bar\A_5$ have a 2-disk suborbifold $D^2/\bar \Z_5$  which is a
triangle: two of its sides are reflection axes meeting in a dihedral point
$\D_5$, the third side is a nonsingular arc (except for its endpoints) which
constitutes the orbifold boundary of $D^2/\bar \Z_5$. The handlebody
orbifold
$\H_5 = \H(\bar\D_5*_{\bar\Z_5}{\bar\A_5})$ is then obtained by
glueing $D^3/\bar\D_5$ and $D^3/\bar\A_5$ along the 2-suborbifolds 
 $D^2/\bar \Z_5$ in their boundaries (or better, by connecting the two 3-disk
orbifolds by a 1-handle orbifold $(D^2/\bar \Z_5) \times [-1,1]$);  note
that the orbifold boundary of this handlebody orbifold is a square 
${\cal D}^2([2,2,2,3])$  whose sides are reflection axes meeting in three dihedral
points $\D_2$ and one $\D_3$.

\medskip

Applying the construction to the first four amalgamated free products, one
obtains the four handlebody orbifolds
$$\H_2 = \H(\bar\D_2*_{\bar\Z_2}{\bar\D_3}), \;\;\;  
\H_3 = \H(\bar\D_3*_{\bar\Z_3}{\bar\A_4}), \;\;\;
\H_4 = \H(\bar\D_4*_{\bar\Z_4}{\bar\S_4}),$$
$$\H_5 = \H(\bar\D_5*_{\bar\Z_5}{\bar\A_5}),$$
and each of these four handlebody orbifolds has the square
${\cal \D}^2([2,2,2,3])$ as its orbifold boundary.

\medskip

For the remaining four free products with amalgamation one constructs in a similar
way the handlebody orbifolds
$$\tilde \H_2 = \H(\D_{2*2}*_{\bar\Z_2}{\bar\D_3}), \;\;\;
\tilde \H_3 = \H(\D_{2*3}*_{\bar\Z_3}{\bar\A_4}), \;\;\;
\tilde \H_4 = \H(\D_{2*4}*_{\bar\Z_4}{\bar\S_4}),$$
$$\tilde \H_5 = \H(\D_{2*5}*_{\bar\Z_5}{\bar\A_5}).$$

For example,  the quotient orbifold $D^3/\D_{2*5}$ is a cone over
its orbifold boundary $S^2/\D_{2*5}$ which is a disk with a singular point
$\Z_2$ in its interior whose boundary consists of a unique reflection axis
starting and finishing in a singular point $\D_5$. Since 
$S^2/\D_{2*5}$ has again a 2-disk suborbifold $D^2/\bar \Z_5$, one  can
construct the handlebody orbifold $\H(\D_{2*5}*_{\bar\Z_5}{\bar\A_5})$
as before.

\medskip

Note that the orbifold boundary of each of the four handlebody orbifolds
$\tilde \H_n$  is a 2-disk ${\cal D}^2(2,[2,3])$ whose boundary
consists of two reflection axes intersecting in two dihedral points $\D_2$
and $\D_3$, and with a singular point $\Z_2$ in its interior; note that 
${\cal D}^2(2,[2,3])$ is a quotient of the square orbifold  ${\cal D}^2([2,2,2,3])$
by a rotational involution $\tau$.

\bigskip

{\bf Remark.} In the case of orientation-preserving actions on orientable handlebodies
and 3-manifolds the maximal order is $12(g-1)$, and the orbifold fundamental groups of
the minimal orientable handlebody orbifolds, of Euler characteristic $-1/12$, are the
following four amalgamated free products:
$$\D_2*_{\Z_2}\D_3, \;\;\; \D_3*_{\Z_3}{\A_4}, \;\;\;
\D_4*_{\Z_4}{\S_4}, \;\;\;  \D_5*_{\Z_5}{\A_5} $$ 
where  $\D_n = (2,2,n)$, 
$\A_4 =(2,3,3)$, $\S_4 = (2,3,4)$ and $\A_5 =
(2,3,5)$ denote the dihedral, tetrahedral,
octahedral and dodecahedral groups, spherical triangle groups which are the
orientation-preserving subgroups of index two in the corresponding extended 
groups (see [Z7]).

\bigskip

Returning to the $G$-action on the 3-manifold $M = V_g \cup_{\partial}V_g'$,
the quotient orbifold  $M/G = (V_g/G) \cup_{\partial}(V_g'/G)$ is
obtained by identifying the minimal handlebody orbifolds
$V_g/G$ and $V_g'/G$ along their boundaries, and both $V_g/G$ and $V_g'/G$
are of one of the eight minimal types described above. Since the boundary
of an orbifold of type $\H_n$ is not homeomorphic to that of an orbifold
of type $\tilde \H_m$, both orbifolds $V_g/G$ and $V_g'/G$ are of the
same type.

\medskip

Suppose first that  $V_g/G = \H_n$ and $V_g'/G = \H_m$. The boundary of both 
$\H_n$ and $\H_m$ is the square orbifold ${\cal D}^2([2,2,2,3])$; up to
isotopy, this square has exactly two orbifold homeomorphisms which are
the identity map and a reflection in a diagonal of the square (connecting 
the two opposite vertices of type $\D_2$ and $\D_3$ and exchanging the other
two vertices of type $\D_2$).

\medskip

Suppose that the boundaries of the handlebody orbifolds $\H_n$ and $\H_m$ are
identified by the identity map of ${\cal D}^2([2,2,2,3])$. As explained
before, both handlebody orbifolds $\H_n$ and $\H_m$  are constructed by
identifying two 3-disk orbifolds along a 2-disk suborbifold $D^2/\bar \Z_n$
and $D^2/\bar \Z_m$ in their boundaries. Identifying the boundaries of
$\H_n$ and $\H_m$ by the
identity map, these 2-disk suborbifolds $D^2/\bar \Z_n$
and $D^2/\bar \Z_m$ fit together along their
boundaries and create a 2-disk suborbifold of $M/G$ whose boundary
consists of two reflection axes meeting in dihedral points $\bar \Z_n \cong
\D_n$ and $\bar \Z_m \cong \D_m$. If $n \ne m$, this 2-disk is a
bad 2-orbifold, i.e. not covered by a manifold, which is a contradiction since
we have the manifold covering $M$ of $M/G$.  Hence $n = m$ and $M/G$ is
obtained by doubling $\H_n$ along the boundary. Then also $M$ is obtained by
doubling the handlebody $V_g$ and its
$G$-action along the boundary, so we are in part i) of the Theorem.

\medskip

Suppose then that the boundaries of $\H_n$ and $\H_m$
are identified by a reflection in a diagonal of the square ${\cal
D}^2([2,2,2,3])$. The Coxeter group $C(n,m)$  acts on the 3-sphere, Euclidean
or hyperbolic 3-space, and the quotient orbifold of this action is the Coxeter
tetrahedron $\C(n,m)$: the underlying topological space is the 3-disk, and
the singular set consists of the boundary of the tetrahedron (the local group
associated to a point is its stabilizer in the Coxeter group).  The  Coxeter
orbifold $\C(n,m)$ has a 2-suborbifold ${\cal D}^2([2,2,2,3])$ (a square
whose vertices are on the four edges of the tetrahedron with labels 2,2,2 and
3, seperating the two edges with labels  $n$ and $m$), and ${\cal
D}^2([2,2,2,3])$  seperates $\C(n,m)$ into the two handlebody orbifolds $\H_n$
and $\H_m$. So in this case the quotient $M/G$ is the Coxeter tetrahedral orbifold 
$\C(n,m)$ and we are in part ii) of the Theorem.

\bigskip

We are left with the cases $V_g/G = \tilde \H_n$ and $V_g'/G = \tilde \H_m$.
The boundary of each of these minimal  handlebody orbifolds  is the 2-disk
${\cal D}^2(2,[2,3])$, and every orbifold homeomorphism of this 2-orbifold is
isotopic to the identity map or to a reflection in a segment which connects the two
singular points of types
$\D_2$ and $\D_3$ and has the singular point $\Z_2$ in its interior.

\medskip

As in the first case, if the boundaries of $\tilde \H_n$ and $\tilde \H_m$ are
identified by the identity map then either $n \ne m$ and there is a bad
2-suborbifold, or $n=m$ and $M/G$ is obtained by
doubling $\tilde \H_n$ along its boundary and we are in case i) of the
Theorem.

\medskip

Finally suppose that the boundaries of $\tilde \H_n$ and $\tilde \H_m$ are
identified by a reflection. Similar as before, the Coxeter orbifold
$\C(n,m;3,3;2,2)$ has a square 2-suborbifold  
${\cal D}^2([3,3,2,2])$ (seperating the two edges with lables $n$ and $m$)
which is invariant under the involution $\tau$ of $\C(n,m;3,3;2,2)$. 
The projection of ${\cal D}^2([3,3,2,2])$ to the
twisted Coxeter orbifold
$\C_\tau(n,m) = \C(n,m;3,3;2,2)/\tau$ is the quotient ${\cal D}^2([3,3,2,2])/\tau$
which is homeomorphic to the
orbifold ${\cal D}^2(2,[2,3])$ and seperates $\C_\tau(n,m)$ into two
handlebody orbifolds $\tilde \H_n$ and $\tilde \H_m$  (e.g., the quotient of the
1-handle orbifold
$D^2/\bar \Z_5 \times [-1/2,1/2]$  by the involution $\tau$ is the 3-disk orbifold 
$D^3/\D_{2*5}$, and hence the quotient of $\H(\bar\A_5*_{\bar\Z_5}{\bar\A_5})$
by $\tau$ gives the handleboldy orbifold $\tilde \H_5 =
\H(\D_{2*5}*_{\bar\Z_5}{\bar\A_5})$).  So by
identifying $\tilde \H_n$ and
$\tilde \H_m$ along their boundaries we obtain 
the twisted Coxeter orbifold $\C_\tau(n,m)$, and we are in case ii) of the Theorem.

\medskip

This completes the proof of the Theorem.

\bigskip

{\bf 3.  Examples and comments}

\medskip

The maximal possible order of a  $G$-action of a closed 3-manifold $M$ which
leaves invariant a Heegaard surface of genus $g>1$ is $48(g-1)$; in this
maximal case, some element of $G$ has to interchange the two handlebodies of
the Heegaard splitting, and the subgroup of index two preserving both
handlebodies gives a strongly maximally symmetric $G$-action on $M$.

\medskip

Suppose that $M$ is irreducible. Then the subgroup of index two preserving
both handlebodies of the Heegaard splitting is
obtained from a Coxeter group $C(n,m)$ or twisted Coxeter group
$C_\tau(n,m)$ as in the Theorem. By the geometrization of finite group
actions in dimension 3, we can assume that the whole group $G$ acts by
isometries; lifting $G$ to the universal covering, we obtain a group of
isometries of $S^3$, $\R^3$ of $\Bbb H^3$ containing the Coxeter group $C(n,m)$
or the twisted Coxeter group $C_\tau(n,m)$ as subgroup of index two, and in
the second case it contains the Coxeter group  $C(n,m;3,3;2,2) \subset
C_\tau(n,m)$ as a subgroup of index 4.  Now any  2-fold or 4-fold extension of
such a Coxeter group is obtained  by adjoining the symmetry
group $\Z_2$  or $\Z_2 \times \Z_2$ of rotations of the Coxeter tetrahedron 
to the Coxeter group (see [M] for the hyperbolic case, the
other cases are similar), and clearly the presence of such symmetries  requires
$n=m$.

\medskip

We denote by $C_\mu(n,n)$ the twisted group generated by 
$C(n,n) = C(n,n;2,2;2,3)$ 
and the involution $\mu$ of the Coxeter tetrahedron $\C(n,n;2,2;2,3)$ 
which exchanges the two edges with label $n$ and inverts the other four
edges, and by $C_{\tau \mu}(n,n)$ the doubly twisted group generated by 
$C(n,n;3,3;2,2)$ and the involutions $\tau$
and $\mu$ of the Coxeter tetrahedron $\C(n,n;3,3;2,2)$. The
possibilities are then the following:

\medskip

spherical:  \hskip 5mm  $C_\mu(2,2), \;\; C_\mu(3,3);  \;\;\;\; 
C_{\tau \mu}(2,2);$

Euclidean:    \hskip 3mm $C_\mu(4,4); \hskip 10mm   C_{\tau \mu}(3,3)$;

hyperbolic:  \hskip 2mm $C_\mu(5,5); \;\;\;\;\;  C_{\tau \mu}(4,4), \;\;
C_{\tau \mu}(5,5).$

\bigskip

Summarizing, we have:

\bigskip

{\bf Proposition.}   {\sl  Let $M$ be closed, irreducible 3-manifold with
a $G$-action of maximal possible order $48(g-1)$
which leaves invariant the Heegaard surface of a Heegaard splitting of genus
$g>1$ of $M$. Then $M$ is obtained as the quotient of the 3-sphere, Euclidean
or hyperbolic 3-space by a normal subgroup of finite index, acting freely, of
a twisted Coxeter group
$C_\mu(n,n)$ or a doubly twisted Coxeter group $C_{\tau \mu}(n,n)$, and the
$G$-action is obtained as the projection of 
$C_\mu(n,n)$ or $C_{\tau \mu}(n,n)$ to $M$. If $M$ is the 3-sphere, the
possible genera are $g$ = 2, 4 and 6.}

\bigskip

Finally, we discuss an explicit example of a hyperbolic, strongly maximally
symmetric 3-manifold $M$ for which, moreover, also the bound  $48(g-1)$ is
obtained: this is the hyperbolic 3-manifold $M$ considered in [Z5], [Z7], 
see also  [Z6] for some further properties.  The universal covering group of
$M$ is a normal torsionfree subgroup of smallest possible index 120 in the
Coxeter group  $C(5,5;3,3;2,2)$, and of index $240 = 24(g-1)$ in its twisted
version $C_\tau(5,5)$.  In particular, $M$ is a strongly maximally symmetric
$G$-manifold of genus $g=11$ (and, by [Z5], also 
the ordinary Heegaard genus of $M$ is equal to 11). Moreover, it follows
from [Z6, Propositions 3.2 and 3.3] that the universal covering group of $M$
is a normal subgroup also of the doubly twisted Coxeter group  $C_{\tau \mu}(n,n)$,
so $C_{\tau \mu}(n,n)$ projects to an isometry group 
of order $480 = 48(g-1)$ of $M$ (which is, in fact, the full isometry group
of $M$), and we are in the situation of the Proposition.

\medskip

We believe that the genus
$11$ of $M$ is the smallest equivariant genus of a hyperbolic $G$-manifold for
which the bound $48(g-1)$  in the  Proposition is  obtained and, more
generally, also the smallest genus of any strongly maximally symmetric
hyperbolic 3-manifold.  For this, one has to check the minimal indices of the
torsionfree normal subgroups of the other hyperbolic Coxeter and twisted
Coxeter groups.

\bigskip \bigskip

\centerline {\bf References}

\medskip

\item {[CM]}  H.S.M. Coxeter, W.O.J. Moser, {\it  Generators and Relations
for Discrete Groups.}  Ergebnisse der Mathematik und ihrer Grenzgebiete 14,
Springer-Verlag 1984

\item {[M]} A.D. Mednykh,  {\it Automorphism groups of three-dimensional
hyperbolic manifolds.}  Soviet Math. Dokl. 32  (1985), 633-636

\item {[MMZ]} D. McCullough, A. Miller, B. Zimmermann,  {\it Group actions on
handlebodies.}  Proc. London Math. Soc.  59   (1989), 373-415

\item {[T1]} W. Thurston, {\it The Geometry and Topology of 3-Manifolds.} 
Lecture Notes Princeton Univ. 1977

\item {[T2]} W. Thurston, {\it Three-dimensional Geometry and Topology.}
Revised and extended version of a part of [T1], Princeton 1990

\item {[WWZ]}  C. Wang, S. Wang, Y. Zhang, {\it Maximum orders of
extendable actions on surfaces.}  Acta Math. Sin. (Engl. Ser.) 32 (2016), 
54-68

\item {[WZ1]}  C. Wang, S. Wang, Y. Zhang, B. Zimmermann, {\it Extending
finite group actions on surfaces over $S^3$.}  Top. Appl. 160 (2013),
2088-2103

\item {[WZ2]}  C. Wang, S. Wang, Y. Zhang, B. Zimmermann, {\it Embedding
surfaces into $S^3$ with maximum symmetry.} 
Groups, Geometry, and Dynamics 9  (2015),  1001-1045

\item {[Z1]} B. Zimmermann,  {\it \"Uber Abbildungsklassen von
Henkelk\"orpern.}  Arch. Math. 33  (1979),  379-382

\item {[Z2]}  B. Zimmermann, {\it Finite group actions on handlebodies and
equivariant Heegaard genus for 3-manifold.}  Top. Appl. 43 (1992),
263-274

\item {[Z3]}  B. Zimmermann, {\it  Finite maximal orientation reversing
group actions on handlebodies and 3-manifolds.}  Rend. Circ. Mat. Palermo 48 
(1999),  549-562

\item {[Z4]} B. Zimmermann,  {\it A survey on large finite group actions on
graphs, surfaces and 3-manifolds.}  Rend. Circ. Mat. Palermo 52  (2003), 47-56

\item {[Z5]} B. Zimmermann,  {\it On a hyperbolic 3-manifold with some
special properties.}  Math. Proc. Camb. Phil. Soc. 113 (1993), 87-90

\item {[Z6]} B. Zimmermann,  {\it On 
maximally symmetric hyperbolic 3-manifolds.}  
Progress in Knot Theory and Related Topics, Travaux En Cours, Hermann, Paris
1997,  143-153

\item {[Z7]} B. Zimmermann,  {\it Genus actions of finite groups on 3-manifolds.} Mich.
Math. J. 43 (1996),  593-610

\bye